\documentclass[reqno,11pt]{amsart}

\newtheorem{theorem}{Theorem}[section]
\newtheorem{lemma}[theorem]{Lemma}
\newtheorem{corollary}[theorem]{Corollary}

\theoremstyle{definition}
\newtheorem{definition}[theorem]{Definition}

 \theoremstyle{remark}
\newtheorem{remark}[theorem]{Remark}

\newcommand\bL{\mathbb{L}}

\newcommand\bR{\mathbb{R}}

\newcommand\bW{\mathbb{W}}

\newcommand\cB{\mathcal{B}}

\newcommand\cD{\mathcal{D}}
\newcommand\cF{\mathcal{F}}

\newcommand\cP{\mathcal{P}}

\newcommand\cU{\mathcal{U}}
\newcommand\cW{\mathcal{W}}

\newcommand{\mysection}[1]{\section{#1}
\setcounter{equation}{0}}

\newcommand\cbrk{\text{$]$\kern-.15em$]$}} 
\newcommand\opar{
\text{\,\raise.2ex\hbox{${\scriptstyle |}$}\kern-.34em$($}} 
\newcommand\cpar{%
\text{$)$\kern-.34em\raise.2ex\hbox{${\scriptstyle |}$}}\,}
\newcommand\obrk{\text{$[$\kern-.15em$[$}}

\begin{document}

\title[It\^o's formula]
{It\^o's formula
for the $L_{p}$-norm of
  stochastic $W^{1}_{p}$-valued processes}
\author{N.V. Krylov}
\address{127 Vincent Hall, University of Minnesota, Minneapolis,
 MN, 55455}
\thanks{The work   was partially supported by
NSF Grant DMS-0653121}
\email{krylov@math.umn.edu}

\renewcommand{\subjclassname}{%
\textup{2000} Mathematics Subject Classification}


\begin{abstract}
We prove It\^o's formula
for the $L_{p}$-norm of
a stochastic $W^{1}_{p}$-valued processes
appearing in the theory of SPDEs in
divergence form.

\end{abstract}

\maketitle

\mysection{Introduction}

Let $(\Omega,\cF,P)$ be a complete probability space
with an increasing filtration $\{\cF_{t},t\geq0\}$
of complete with respect to $(\cF,P)$ $\sigma$-fields
$\cF_{t}\subset\cF$. Denote by $\cP$ the predictable
$\sigma$-field in $\Omega\times(0,\infty)$
associated with $\{\cF_{t}\}$. Let
 $w^{k}_{t}$, $k=1,2,...$, be independent one-dimensional
Wiener processes with respect to $\{\cF_{t}\}$.  
Let $\cD$ be the space 
of
generalized functions on the Euclidean 
$d$-dimensional space $\bR^{d}$
of points $x=(x^{1},...,x^{d})$. 
We consider processes with values in $\cD$ whose
stochastic differential is given by 
\begin{equation}
                                       \label{12.3.1}
du_{t}=(D_{i}f^{i}_{t}+f^{0}_{t})\,dt
+g^{k}_{t}\,dw^{k}_{t},
\end{equation}
where $f^{j}_{t},g^{k}_{t}$ are $L_{p}$-valued processes, $u_{t}$
is a $W^{1}_{p}$-valued process,
and the summation convention over repeated indices is enforced.
Our main goals are to give  conditions
on $u,f^{j}$, and $g^{k}$, which are sufficient
to assert that $u_{t}$ is a continuous $L_{p}$-valued 
process, and to derive  It\^o's formula for
$\|u_{t}\|_{L_{p}}^{p}$.

This was never done before, no matter how strange
it may look. The hardest step   is 
showing that $u_{t}$ is continuous as an $L_{p}$-valued
function.  More or less standard fact is that 
under natural conditions
one can estimate
\begin{equation}
                                                     \label{5.3.3}
E\sup_{t}\|u_{t}\|_{L_{p}}^{p}
\end{equation}
and from here and
equation \eqref{12.3.1}, implying that $(u_{t},\varphi)$
is continuous in $t$ for any test function
$\varphi\in C^{\infty}_{0}$, one used to derive that $u_{t}$
is only a weakly continuous $L_{p}$-valued process.
Even though the above mentioned It\^o's formula
was not proved,
the fact that, actually,
$u_{t}$ is indeed continuous as an $L_{p}$-valued
process was known and proved by different methods
for $p=2$, on the basis of abstract results
for SPDEs in Hilbert spaces, and for $p>2$,
 on the basis of embedding theorems
for stochastic Banach spaces. In this way of arguing 
proving the continuity of $\|u_{t}\|_{L_{p}}^{p}$
 required a full blown theory
of SPDEs with constant coefficients (cf.~\cite{Ro}
and \cite{Kr99}).
We present a ``direct" and self-contained
proof of the formula and the continuity.

Finally, we mention that there are many situations
in which It\^o's formula is known for Banach space valued processes.
See, for instance, \cite{Br} and the references therein.
These formulas could be more general in some respects but they
do not cover our situation and are closer to our Lemma
\ref{lemma 4.17.1} where the term $D_{i}f^{i}$ is not present in
\eqref{12.3.1}.

\mysection{Main result}

We take a stopping time $\tau$ and fix a number 
$$
p\geq2. 
$$
Denote $L_{p}=L_{p}(\bR^{d})$.
We use the same notation $L_{p}$ for vector- and matrix-valued
or else
$\ell_{2}$-valued functions such as
$g_{t}=(g^{k}_{t})$ in \eqref{12.3.1}. For instance,
if $u(x)=(u^{1}(x),u^{2}(x),...)$ is 
an $\ell_{2}$-valued measurable function on $\bR^{d}$, then
$$
\|u\|^{p}_{L_{p}}=\int_{\bR^{d}}|u(x)|_{\ell_{2}}^{p}
\,dx
=\int_{\bR^{d}}\big(
\sum_{k=1}^{\infty}|u^{k}(x)|^{2}\big)^{p/2}
\,dx.
$$

Introduce
$$
D_{i}=\frac{\partial}{\partial x^{i}},\quad i=1,...,d.
$$
By $Du$   we mean the gradient  with respect
to $x$ of a function $u$ on $\bR^{d}$.

As usual,  
$$
W^{1}_{p}=\{u\in L_{p}: Du\in L_{p}\},
\quad
 \|u\|_{W^{1}_{p}}=
\|u\|_{L_{p}}+\|Du\|_{L_{p}}.
$$

If $\tau$ is a stopping time, then
$$
\bL _{p}(\tau):=L_{p}(\opar 0,\tau\cbrk,\cP,
L_{p}),\quad
\bW^{1}_{p}(\tau):=L_{p}(\opar 0,\tau\cbrk,\cP,
W^{1}_{p}).
$$
We also need the space $\cW^{1}_{p}(\tau)$,
which is the space of functions $u_{t}
=u_{t}(\omega,\cdot)$ on $\{(\omega,t):
0\leq t\leq\tau,t<\infty\}$ with values
in the space of generalized functions on $\bR^{d}$
and having the following properties:

(i) We have $u_{0}\in L_{p}(\Omega,\cF_{0},L_{p})$;

(ii)  We have $u
\in \bW^{1}_{p}(\tau )$;

(iii) There exist   $f^{i}\in \bL_{p}(\tau)$,
$i=0,...,d$, and $g=(g^{1},g^{2},...)\in \bL_{p}(\tau)$
such that
 for any $\varphi\in C^{\infty}_{0}$ with probability 1
for all  $t\in[0,\infty)$
we have
$$
(u_{t\wedge\tau},\varphi)=(u_{0},\varphi)
+\sum_{k=1}^{\infty}\int_{0}^{t}I_{s\leq\tau}
(g^{k}_{s},\varphi)\,dw^{k}_{s}
$$
\begin{equation}
                                                 \label{1.2.1}
+\int_{0}^{t}I_{s\leq\tau}\big(
 (f^{0}_{s},\varphi)-(f^{i}_{s},D_{i}\varphi)\big)\,ds.
\end{equation}
In particular, for any $\phi\in C^{\infty}_{0}$, the process
$(u_{t\wedge\tau},\phi)$ is $\cF_{t}$-adapted and continuous.
In case that property (iii) holds, we say that \eqref{12.3.1}
holds for $t\leq\tau$.

The reader can find 
in \cite{Kr99} a discussion of (ii) and (iii),
in particular, the fact that the series in \eqref{1.2.1}
converges uniformly in probability on every finite
subinterval of $[0,\tau]$. This will also be seen
from the proof of Lemma~\ref{lemma 5.5.3}.

Here is our main result.
\begin{theorem}
                                     \label{theorem 12.3.1}
  Let
$u\in\cW^{1}_{p}(\tau)$, 
$f^{j}\in\bL_{p}(\tau)$,
  $g=(g^{k})\in\bL_{p}(\tau)$ and assume that
\eqref{12.3.1} holds
for $t\leq\tau$ in the sense of generalized functions.
Then there is a set $\Omega'\subset\Omega$ of full probability
such that

(i) $u_{t\wedge\tau}I_{\Omega'}$
is a continuous $L_{p}$-valued $\cF_{t}$-adapted function on
$[0,\infty)$;

(ii) for all $t\in[0,\infty)$ and $\omega\in\Omega'$ It\^o's formula holds:
$$
\int_{\bR^{d}}|u_{t\wedge\tau}|^{p}\,dx
=\int_{\bR^{d}}|u_{0}|^{p}\,dx
+p
\int_{0}^{t\wedge\tau }\int_{\bR^{d}}|u _{s}|^{p-2}
u _{s}
g^{k }_{s}\,dx\,dw^{k}_{s}
$$
$$+
\int_{0}^{t\wedge\tau }
\big( \int_{\bR^{d}}\big[p|u_{t}|^{p-2}u_{t}f^{0}_{t}
-p(p-1)|u_{t}|^{p-2}f^{i}_{t}D_{i}u_{t}
$$
\begin{equation}
                                            \label{4.19.5}
+(1/2)p(p-1)|u_{t}|^{p-2}|g_{t}|_{\ell_{2}}^{2}
\big]\,dx\big)\,dt.
\end{equation}

Furthermore,  for any $T\in[0,\infty)$ and
$$
 E\sup_{t\leq\tau\wedge T}
\|u_{t}\|^{p}_{L_{p}}\leq  2E\|u_{0}\|^{p}_{L_{p}}+
NT^{p-1}\|f^{0}\|^{p}_{\bL_{p}(\tau)}
$$
\begin{equation}
                                         \label{4.11.5}
+NT^{(p-2)/2}(\sum_{i=1}^{d}\|f^{i}\|^{p}_{\bL_{p}(\tau)}
+\|g\|^{p}_{\bL_{p}(\tau)}+\|Du\|^{p}_{\bL_{p}(\tau)}) ,
\end{equation}
where $N=N(d,p)$.  

\end{theorem}

We prove Theorem \ref{theorem 12.3.1}
in Section \ref{section 5.13.2} after we prepare
the necessary tools in Sections \ref{section 5.5.1}-\ref{section 5.13.1}.

Here is an ``energy" estimate.

\begin{corollary}
                               \label{corollary 4.19.1}
Under the conditions of Theorem \ref{theorem 12.3.1}
$$
E\int_{\bR^{d}}|u_{0}| ^{p}\,dx+E\int_{0}^{\tau}
\big( \int_{\bR^{d}}\big[p|u_{t}|^{p-2}u_{t}f^{0}_{t}
-p(p-1)|u_{t}|^{p-2}f^{i}_{t}D_{i}u_{t}
$$
\begin{equation}
                                       \label{12.3.2}
+(1/2)p(p-1)|u_{t}|^{p-2}|g_{t}|_{\ell_{2}}^{2}
\big]\,dx\big)\,dt 
\geq EI_{\tau<\infty}\int_{\bR^{d}}|u_{\tau}| ^{p}\,dx.
\end{equation}
Furthermore, if $\tau$ is bounded then there is an equality
instead of inequality in \eqref{12.3.2}.

\end{corollary}

The proof of the corollary is given in Section \ref{section 5.13.2}.

\mysection{Auxiliary results}
                           
                                    \label{section 5.5.1}

We need two well-known results
(see, for instance, Lemma 6.1 and Corollary 6.2 in \cite{Kr08}),
 which we prove
for completeness of presentation.

\begin{lemma}
                                        \label{lemma 1.15.1}
Let $(E,\Sigma,\mu)$ be a measure space,
$r\in[1,\infty)$, 
$u_{n},u\in L_{r}(\mu)$,  and
$u_{n}\to u$ in measure. Finally, let
$$
\|u_{n}\|_{L_{r}(\mu)}\to\|u\|_{L_{r}(\mu)}.
$$
Then
\begin{equation}
                                               \label{1.16.1}
\|u_{n}-u\|_{L_{r}(\mu)}\to0 .
\end{equation}
\end{lemma}

Proof. We have
$$
|\,|u|^{r}-|u_{n}|^{r}|=(|u|^{r}-|u_{n}|^{r})_{+}
-(|u|^{r}-|u_{n}|^{r}).
$$
Upon integrating through this equation and observing that
$(|u|^{r}-|u_{n}|^{r})_{+}\leq|u|^{r}$ we conclude
 by the dominated convergence theorem that
\begin{equation}
                                               \label{1.16.2}
\int_{E}|\,|u|^{r}-|u_{n}|^{r}|\,\mu(dx)\to0.
\end{equation}

Next, if $|u_{n}-u|\geq 3|u|$, then $|u_{n}|+|u|\geq3|u|$,
$|u|\leq(1/2)|u_{n}|$, $|u|^{r}\leq(1/2)|u_{n}|^{r}$,
$$
|u_{n}|^{r}-|u|^{r}\geq(1/2)|u_{n}|^{r},
\quad |u_{n}-u|\leq|u_{n}|+|u|\leq2|u_{n}|,
$$
$$
|u_{n}-u|^{r}\leq2^{r}|u_{n}|^{r}
\leq 4^{r} (|u_{n}|^{r}-|u|^{r}),
$$
which along with \eqref{1.16.2}
imply  that
$$
\int_{E}|u_{n}-u|^{r}I_{|u_{n}-u|\geq 3|u|}\,\mu(dx)\to0.
$$
Furthermore,
$$
\int_{E}|u_{n}-u|^{r}I_{|u_{n}-u|<3|u|}\,\mu(dx)
\to0
$$
by the dominated convergence theorem. By combining the 
last two relations we 
come to \eqref{1.16.1}. The lemma is proved.
\begin{corollary}
                              \label{corollary 1.15.1}
Let $(E,\Sigma,\mu)$ be a measure space,
$r,s\in(1,\infty)$, $r^{-1}+s^{-1}=1$,
$u_{n},u\in L_{r}(\mu)$, $v_{n},v\in L_{s}(\mu)$,
$u_{n}\to u$ and $v_{n}\to v$ in measure. Finally, let
$$
\|u_{n}\|_{L_{r}(\mu)}\to\|u\|_{L_{r}(\mu)},\quad
\|v_{n}\|_{L_{s}(\mu)}\to\|u\|_{L_{s}(\mu)}.
$$
Then
$$
\int_{E}|u_{n}v_{n}-uv|\,\mu(dx)\to0,\quad
\int_{E}u_{n}v_{n}\,\mu(dx)\to\int_{E}uv\,\mu(dx).
$$
\end{corollary}

Indeed, it suffices to use H\"older's inequality and  
the formula
$$
u_{n}v_{n}-uv=(u_{n}-u)v+(v_{n}-v)u+(u_{n}-u)(v_{n}-v).
$$

\mysection{Integrating $\bL_{p}$ functions}
                           
                                    \label{section 5.3.2}

Most likely a big part of what follows in this section
can be obtained from some abstract constructions in \cite{Br}.
However, it does not look easy to obtain estimate
\eqref{5.9.1}. In any case, it is worth giving
all rather simple arguments for completeness. Set
$$
\bL_{p}=\bL_{p}(\infty)
$$
and for Borel subsets $\Gamma$ of a Euclidean space denote
by $\cB(\Gamma)$ the $\sigma$-field of Borel subsets of $\Gamma$.

\begin{definition}
                                      \label{definition 5.9.1}
By $\cU_{p}$ we denote the set of functions
$u=u_{t}(x)=u_{t}(\omega,x)$ on $\Omega\times[0,\infty)
\times\bR^{d}$ such that

(i) $u$ is measurable with respect to $\cF\otimes\cB([0,\infty))
\otimes\cB(\bR^{d})$;

(ii) for each $x$, the function $u_{t}(x)$ is $\cF_{t}$-adapted;

(iii) $u_{t}(x)$ is continuous  in $t\in[0,\infty)$ for each $(\omega,x)$;
 
(iv) the function $u_{t}(\omega,\cdot)$
 as a function of $(\omega,t)$ is $L_{p}$-valued,
$\cF_{t}$-adapted, and continuous in $t$ for any $ \omega$.
 
\end{definition}

\begin{lemma}
                                     \label{lemma 5.5.3}
 Let
  $g=(g^{k})\in\bL_{p} $.
Then there exists a function $u\in \cU_{p}$
such that
for any $\phi\in C^{\infty}_{0}$ the equation
\begin{equation}
                                                  \label{5.9.2}
(u_{t},\phi)=\sum_{k=1}^{\infty}\int_{0}^{t}(g^{k}_{s},
\phi)\,dw^{k}_{s}
\end{equation}
holds for all $t\in[0,\infty)$ with probability one.
Furthermore, for any $T\in[0,\infty)$ we have
\begin{equation}
                                                  \label{5.9.1}
E\int_{\bR^{d}}\sup_{t\leq T}
|u _{t }(x)|^{p}\,dx\leq NT^{(p-2)/2}E\int_{0}^{T }
\|g _{s}\| ^{p }_{L_{p}}\,ds,
\end{equation}
where $N=N(p)$.

\end{lemma}

Proof. First assume that there is an integer $j\geq1$,
 (nonrandom) functions $g^{ik}\in C^{\infty}_{0}$, 
and bounded stopping times $\tau_{0}\leq
\tau_{1}\leq...\leq\tau_{j}$
such that $g^{k}\equiv0$ for $k>j$ and
$$
g^{k}_{t}(x)=\sum_{i=1}^{j}
g^{ik}(x)I_{(\tau_{i-1},\tau_{i}]}(t)
$$
for $k\leq j$.

Then define
$$
u_{t}(x)=\sum_{i,k=1}^{j}g^{ik}(x)(w^{k}_{t\wedge\tau_{i}}-
w^{k}_{t\wedge\tau_{i-1}}).
$$
Obviously, $u\in\cU_{p}$. Furthermore, \eqref{5.9.2} holds
for any $\phi\in C^{\infty}_{0}$ for all $t$ with probability one
since its right-hand side equals
$$
\sum_{i,k=1}^{j}(g^{ik},\phi)(w^{k}_{t\wedge\tau_{i}}-
w^{k}_{t\wedge\tau_{i-1}})
$$
for all $t$ with probability one. Next, by the Burkholder-Davis-Gundy
inequalities for each $x$
$$
E\sup_{t\leq T}
|u _{t }(x)|^{p} =E\sup_{t\leq T}\big|\sum_{k}\int_{0}^{t}
g^{k}_{s}(x)\,dw^{k}_{s}\big|^{p}\leq N
E\big(\int_{0}^{T}|g_{s}(x)|_{\ell_{2}}^{2}\,ds\big)^{p/2},
$$
which after applying H\"older's inequality ($p\geq2$) yields
$$
E\sup_{t\leq T}
|u _{t }(x)|^{p}  \leq NT^{(p-2)/2}
E \int_{0}^{T}|g_{s}(x)|_{\ell_{2}}^{p}\,ds.
$$
We integrate this inequality over $\bR^{d}$ and use the fact that
the measurability properties of $g,u$ and the continuity
of $u_{t}$ in $t$ allow us to use Fubini's theorem.
Then we come to \eqref{5.9.1}.

By Theorem 3.10 of \cite{Kr99} the set of $g$'s like the one above
is dense in $\bL_{p}$. Therefore, to prove the lemma
it suffices to show that the set of $g$'s for which
the statements of the lemma are true is closed in $\bL_{p}$.

Take a sequence $g^{n}=(g^{nk} )\in\bL_{p}$, $n=1,2,...$,
such that for each $n$ there is a function $u^{n}$
corresponding to $g^{n}$ and possessing the asserted properties.
Assume that for a $g\in\bL_{p}$ we have
$g^{n}\to g$ in $\bL_{p}$ as $n\to\infty$.
Using a subsequence of $g^{n}$ we may assume that for any $T\in
[0,\infty)$
\begin{equation}
                                                  \label{5.9.4}
E\int_{\bR^{d}}\sup_{t\leq T}
|u^{n+1}_{ t }(x)-u _{t}^{n}(x)|^{p}\,dx\leq  T^{(p-2)/2}2^{-n}.
\end{equation}
Introduce  
$$
A_{n}=\{(\omega,x):
\sup_{t\leq n}
|u^{n+1}_{ t }(x)-u _{t}^{n}(x)|\geq n^{-2}\}.
$$
Then  
$$
\sum_{n=1}^{\infty}
E\int_{\bR^{d}}\sup_{t\leq n}
|u^{n+1}_{ t }(x)-u _{t}^{n}(x)|I_{A_{n}}(x)\,dx
$$
$$
\leq 
\sum_{n=1}^{\infty}n^{2(p-1)}E\int_{\bR^{d}}\sup_{t\leq n}
|u^{n+1}_{ t }(x)-u _{t}^{n}(x)|^{p}\,dx
$$
$$
\leq  \sum_{n=1}^{\infty}
n^{2(p-1)}n^{(p-2)/2}2^{-n}<\infty,
$$
implying that
$$
\sum_{n=1}^{\infty}\sup_{t\leq n}
|u^{n+1}_{ t }(x)-u _{t}^{n}(x)|I_{A_{n}}(x)<\infty
$$
for almost all $(\omega,x)$. The series with the complements
of $A_{n}$ in place of $A_{n}$ obviously converges
everywhere. We conclude that the $\cF\otimes\cB(\bR^{d})$-measurable set  
$$
G=\{(\omega,x):\sum_{n=1}^{\infty}\sup_{t\leq n}
|u^{n+1}_{ t }(x)-u _{t}^{n}(x)|<\infty\}
$$
has full measure. By Fubini's theorem the function
$P((\omega,x)\in G)$ is a Borel function of $x$
equal to 1 for almost all $x$. Accordingly we introduce
a Borel set of full measure
$$
\Gamma=\{x:P((\omega,x)\in G)=1\}
$$
and the $\cF\otimes\cB(\bR^{d})$-measurable set  $G'$ of full measure by
$$
G'=\{(\omega,x):x\in\Gamma,\sum_{n=1}^{\infty}\sup_{t\leq n}
|u^{n+1}_{ t }(x)-u _{t}^{n}(x)|<\infty\}.
$$

Now define
\begin{equation}
                                                  \label{5.9.3}
u'_{t}(x)=\lim_{n\to\infty}u^{n}_{t}(x)
\end{equation}
for   $(\omega, x)\in G'$, $t\geq0$ 
and set $u'_{t}(x)\equiv0$ for   $(\omega, x)\not\in G'$.
Also set
$$
v_{t}(x)=\begin{cases}
\lim_{n\to\infty}u^{n}_{t}(x)&\text{if the limit exists},
\\
0&\text{otherwise}.
\end{cases}
$$
Obviously,
\begin{equation}
                                                  \label{6.6.1}
u'_{t}(\omega,x)=v_{t}(\omega,x)I_{G'}(\omega,x).
\end{equation}
Furthermore, $v$ is known to be $\cF\otimes\cB([0,\infty))
\otimes\cB(\bR^{d})$-measurable since $u^{n}$ possess this property.
It follows that $u'$ is $\cF\otimes\cB([0,\infty))
\otimes\cB(\bR^{d})$-measurable. For each $x$, the functions
$u^{n}_{t}(x)$ are $\cF_{t}$-adapted and so is $v_{t}(x)$. Also
$I_{G'}(\omega,x)$ is $\cF_{0}$-measurable (and hence
$\cF_{t}$-adapted) for each $x$ since the
$\cF_{t}$ are complete and
$$
P(I_{G'}(\omega,x)=1)=P((\omega,x)\in G,x\in\Gamma)
$$
equals zero if $x\not\in\Gamma$ and one if $x\in\Gamma$
by the choice of $\Gamma$. Now equation \eqref{6.6.1} allows us to conclude that
$u'_{t}(x)$ is $\cF_{t}$-adapted for each $x$.

Since the limit in \eqref{5.9.3} is uniform in $t$
on any finite interval, we see that $u'_{t}$ is continuous
in $t$ for any $(\omega,x)$. In particular,
$$
\sup_{t\leq T}|u'_{t}(x)|
$$
is $\cF\otimes\cB(\bR^{d})$-measurable, estimate \eqref{5.9.1}
with $u'$ in place of $u$
makes sense and holds owing to Fatou's lemma and the assumption
on $u^{n}$.

Estimate \eqref{5.9.1} shows that there is a set $\Omega'\in\cF_{0}$
of full probability such that $u'_{t}(\omega,\cdot)
\in L_{p}$ for all $t$
if $\omega\in\Omega'$ and moreover
$$
\int_{\bR^{d}}\sup_{t\leq T}
|u' _{t }(\omega,x)|^{p}\,dx<\infty
$$
for any $T\in[0,\infty)$ if $\omega\in\Omega'$.
This fact, the continuity of $u'_{t}$ in $t$,
and the dominated convergence theorem imply that
$u'_{t}$
is continuous as an $L_{p}$-valued function of $t$
for any $ \omega\in\Omega'$. We now set
$$
u_{t}(\omega,x)=u'_{t}(\omega,x)I_{\Omega'}(\omega).
$$
Then we see that to show that $u\in\cU_{p}$
it suffices to prove that $u_{t}(\omega,\cdot)$
is $\cF_{t}$-adapted as an $L_{p}$-valued function.

Obviously, to do this step it suffices to
 prove the assertion of the lemma related to \eqref{5.9.2}.
 By the Burkholder-Davis-Gundy
inequalities for any $T\in[0,\infty)$ 
$$
 E\sup_{t\leq T}\big|\sum_{k}\int_{0}^{t}
 (g^{k}_{s}-g^{nk}_{s},\phi) \,dw^{k}_{s}\big|^{p}\leq N
E\big(\int_{0}^{T}
\sum_{k=1}^{\infty}(g^{k}_{s}-g^{nk}_{s}, \phi )^{2}\,ds\big)^{p/2}
$$
$$
\leq N
E\big(\int_{0}^{T}
\sum_{k=1}^{\infty}(|g^{k}_{s}-g^{nk}_{s}|^{2},|\phi|)
\|\phi\|^{2}_{L_{2}}\,ds\big)^{p/2}
$$
$$
\leq N
E\big(\int_{0}^{T}\int_{\bR^{d}}
  |g _{s}-g^{n }_{s}|_{\ell_{2}}^{2} |\phi| 
 \,dxds\big)^{p/2}
$$
\begin{equation}
                                                      \label{5.12.2}
\leq N
E \int_{0}^{T}\int_{\bR^{d}}
  |g _{s}-g^{n }_{s}|_{\ell_{2}}^{p} \,dx
 \,ds \leq N\|g-g^{n}\|_{\bL_{p}},
\end{equation}
where $N$ is independent of $n$. In addition, estimate \eqref{5.9.4}
easily imply that
$$
E\sup_{t\leq T}|(u_{t}-u^{n}_{t},\phi)|\to0
$$
as $n\to\infty$ for any $T\in[0,\infty)$. By combining
these fact and passing to the limit in
\eqref{5.9.2} with $u^{n}$ in place of $u$ we get the desired result
and the lemma is proved. 

\begin{remark}
It is tempting to assert that $u_{t}(x)$ is $\cP\otimes\cB(\bR^{d})$-measurable
since it is $\cF\otimes\cB([0,\infty))
\otimes\cB(\bR^{d})$-measurable and,
for each $x$, it is predictable. However, we do not know
if this assertion is true.

\end{remark}

In a similar way, without using the Burkholder-Davis-Gundy
inequalities, the following result is established.

\begin{lemma}
                                     \label{lemma 5.9.1}
 Let
  $f\in\bL_{p}(\infty)$.
Then there exists a function $u\in \cU_{p}$
such that
for any $\phi\in C^{\infty}_{0}$ the equation
\begin{equation}
                                                  \label{5.9.5}
(u_{t},\phi)= \int_{0}^{t}(f_{s},
\phi)\,ds
\end{equation}
holds for all $t\in[0,\infty)$ with probability one.
Furthermore, for any $T\in[0,\infty)$ we have
\begin{equation}
                                                  \label{5.9.6}
E\int_{\bR^{d}}\sup_{t\leq T}
|u _{t }(x)|^{p}\,dx\leq NT^{p-1}E\int_{0}^{T }
\|f_{s}\| ^{p }_{L_{p}}\,ds,
\end{equation}
where $N=N(p)$.

\end{lemma}
\begin{remark}
Observe that the integral on the right in \eqref{5.9.5}
need not exist for each $\omega$ since $(f_{t},\phi)$
is generally only measurable with respect to
the completion of $\cP$ and 
the function $(f_{t},\phi)$ need not be Lebesgue measurable 
in $t$ for each  $\omega$. In case of 
Lemma \ref{lemma 5.5.3} this moment
does not arise because of freedom in defining the stochastic
integrals.
\end{remark}

\mysection{It\^o's formula in a simple case}
                                              \label{section 5.13.1}

The goal of this section is to prove the following result.     

\begin{lemma}
                                     \label{lemma 4.17.1}
  Let
$f \in\bL_{p}(\tau)$,
  $g=(g^{k})\in\bL_{p}(\tau)$ and assume that
we are given a function $u_{t}$ on $\Omega\times[0,\infty)$ with values
in the space of distributions on $\bR^{d}$ such that
$u_{0}\in L_{p}(\Omega,\cF_{0},L_{p})$   and
for any $\phi\in C^{\infty}_{0}$ with probability one
for all $t\in[0,\infty)$ we have
\begin{equation}
                                       \label{4.18.6}
(u_{t\wedge\tau},\phi)=(u_{0},\phi)
+\int_{0}^{t}I_{s\leq\tau}(f_{s},\phi)\,ds+\sum_{k=1}^{\infty}
\int_{0}^{t}(g^{k}_{s},\phi)I_{s\leq\tau}\,dw^{k}_{s}.
\end{equation}
Then, there is a set $\Omega'\in\cF_{0}$ of full  probability such that 

(i) $u_{t\wedge\tau}I_{\Omega'}$ is an $L_{p}$-valued $\cF_{t}$-adapted
continuous process on $[0,\infty)$,

(ii) for all $t\in[0,\infty)$ and $\omega\in\Omega'$
$$
\|u_{t\wedge\tau}\|^{p}_{L_{p}}
= \|u_{0}\|^{p}_{L_{p}} + \int_{0}^{t\wedge\tau }
 \big[p\int_{\bR^{d}}|u _{s}|^{p-2}
u _{s}f _{s}\,dx\,ds
$$
$$
+(1/2)p(p-1)\int_{\bR^{d}}|u _{s}|^{p-2}
|g _{s}|_{\ell_{2}}^{2}\,dx
\big]\,ds
$$
\begin{equation}
                                       \label{4.18.5}
+p
\int_{0}^{t\wedge\tau }\int_{\bR^{d}}|u _{s}|^{p-2}
u _{s}
g^{k }_{s}\,dx\,dw^{k}_{s}.
\end{equation}

\end{lemma}

Proof. First observe that
the right-hand sides of \eqref{4.18.6}  and \eqref{4.18.5}
will be affected only on the set of probability zero independent 
of $t$ if we replace
 $f$ and $g$ with
$L_{p}$-valued {\em predictable\/} functions $\hat{f}$ and $\hat{g}$
such that 
$$
\|f_{t}-\hat{f}_{t}\|_{L_{p}}+\|g_{t}-\hat{g}_{t}\|_{L_{p}}=0
$$
for almost all $(\omega,t)$. It follows that
without losing generality we may assume that
 $f$ and $g$ are predictable as $L_{p}$-valued functions.

 Lemmas \ref{lemma 5.5.3} and \ref{lemma 5.9.1}
allow us to find a $v\in\cU_{p}$ such that for any 
$\phi\in C^{\infty}_{0}$ equation
\eqref{4.18.6} with $v_{t}$ in place of
$u_{t\wedge\tau}$ holds for all $t$ with probability one.
It follows that for any countable set $A\subset C^{\infty}_{0}$
there exists a set $\Omega'$ of full probability
such that for any $\omega\in\Omega'$, $\phi\in A$, and $t\geq0$
we have $(u_{t\wedge\tau},\phi)=(v_{t},\phi)$. If the set $A$
is chosen appropriately, then we conclude that 
$u_{t\wedge\tau}=v_{t}$ in the sense
the distributions, whenever $\omega\in\Omega'$ and $t\geq0$.
 In particular, assertion (i) holds with this $\Omega'$. 

This argument allows us to assume that $\tau=\infty$ and
 $u\in\cU_{p}$
and concentrate on proving \eqref{4.18.5}. This argument
also shows that
for any $T\in[0,\infty)$
$$
E\int_{\bR^{d}}\sup_{t\leq T}
|u_{t }(x)|^{p}\,dx<\infty
$$
implying that there exists a set $\Omega''$ of full probability
such that for any $\omega\in\Omega''$ we have
\begin{equation}
                                                     \label{5.9.7}
 \int_{\bR^{d}}\sup_{t\leq T}
|u_{t }(x)|^{p}\,dx<\infty,\quad\forall T\in[0,\infty).
\end{equation}

Now, take a nonnegative
 function $\zeta\in C^{\infty}_{0}(\bR^{d})$ with unit integral, 
for $\varepsilon>0$,
define $\zeta_{\varepsilon}=\varepsilon^{-d}
\zeta(x/\varepsilon)$,
and for any locally summable $h$ given on $\bR^{d}$
introduce the notation
 $
h^{(\varepsilon)}=h*\zeta_{\varepsilon}$.
Then \eqref{4.18.6} implies that for each $x$
almost surely for all  $t\in[0,\infty)$
\begin{equation}
                                                     \label{4.16.1}
 u^{(\varepsilon)}_{t}(x)
= u^{(\varepsilon)}_{0}(x)+
\int_{0}^{t}  f^{(\varepsilon)}_{s}(x)\,ds
+ \int_{0}^{t}g^{k(\varepsilon)}_{s}(x) \,dw^{k}_{s}.
\end{equation}
By It\^o's formula, for each $x$ 
$$
|u^{(\varepsilon)}_{t }|^{p}
=|u^{(\varepsilon)}_{0}|^{p}+
\int_{0}^{t }p|u^{(\varepsilon)}_{s}|^{p-2}
u^{(\varepsilon)}_{s}
g^{k(\varepsilon)}_{s}\,dw^{k}_{s} 
$$
\begin{equation}
                                                     \label{4.12.1}
+\int_{0}^{t }
 \big[p|u^{(\varepsilon)}_{s}|^{p-2}
u^{(\varepsilon)}_{s}f^{  (\varepsilon) }_{s}\,ds
+(1/2)p(p-1)|u^{(\varepsilon)}_{s}|^{p-2}
|g^{(\varepsilon)}_{s}|_{\ell_{2}}^{2}
\big]\,ds
\end{equation}
(a.s.),
where we dropped the argument $x$ for simplicity.
We want to integrate this equality over $\bR^{d}$
and use the stochastic and deterministic Fubini's
theorems. We will see that there is no difficulties
with the integral with respect to $ds$.
However, in order to be able to apply
the stochastic version of Fubini's theorem
we need at least that the resulting stochastic integral
make sense, that is we need at least the inequality
$$
\int_{0}^{t}\sum_{k=1}^{\infty}
\big(\int_{\bR^{d}}|u^{(\varepsilon)}_{s}|^{ p-1}
|g^{k (\varepsilon)}_{s}|  \,dx\big)^{2}\,ds<\infty
$$
to hold (a.s.). The computations below show that, actually,
for a sequence of stopping times $\tau_{n}\uparrow\infty$,
\begin{equation}
                                                     \label{4.16.3}
E\int_{0}^{t\wedge\tau_{n}}\sum_{k=1}^{\infty}
\big(\int_{\bR^{d}}|u^{(\varepsilon)}_{s}|^{ p-1}
|g^{k (\varepsilon)}_{s}|  \,dx\big)^{2}\,ds<\infty
\end{equation}
and this is known to be sufficient  to apply
the stochastic version of Fubini's theorem.
By the way, notice that $u^{(\varepsilon)}_{t}(x)$
is continuous (infinitely differentiable) in $x$ for 
any $(\omega,t)$. Therefore, it is $\cF_{t}\otimes\cB(\bR^{d})$-measurable. 
Since it is also continuous in $t$ for each $(\omega,x)$,
the function $u^{(\varepsilon)}_{t}(x)$ if 
$\cP\otimes\cB(\bR^{d})$-measurable
and there is no measurability obstructions in applying
Fubini's theorems.

 To deal with the integral with respect to $s$ observe that
by Young's inequality for any $t\in[0,\infty)$
$$
\int_{0}^{t}
 |u^{(\varepsilon)}_{s}|^{p-1}
|f^{(\varepsilon)}_{s}|\,ds\leq\frac{\gamma^{p/(p-1)}}{t}
\int_{0}^{t }|u^{(\varepsilon)}_{s}|^{p }\,ds
+\frac{t^{p-1}}{\gamma^{p}}
\int_{0}^{t }|f^{(\varepsilon)}_{s}|^{p }\,ds,
$$
\begin{equation}
                                                   \label{4.18.2}
\int_{0}^{t}|u^{(\varepsilon)}_{s}|^{p-2}
|g^{(\varepsilon)}_{s}|_{\ell_{2}}^{2}\,ds
\leq \frac{\gamma^{p/(p-2)}}{t}
\int_{0}^{t}|u^{(\varepsilon)}_{s}|^{p }\,ds
+\frac{t^{(p-2)/2}}{\gamma^{p/2}}
\int_{0}^{t}|g^{(\varepsilon)}_{s}|_{\ell_{2}}^{p }\,ds,
\end{equation}
where  $\gamma>0$ is any number (however, if $p=2$
we set $\gamma=1$ in the second inequality). Actually, below in this proof
we only need \eqref{4.18.2} with $\gamma=1$.
More general $\gamma$'s will appear in the proof of Theorem
\ref{theorem 12.3.1}.

 By Minkowski's inequality
\begin{equation}
                                            \label{4.22.1}
\sum_{k=1}^{\infty}
\big(\int_{\bR^{d}}|u^{(\varepsilon)}_{s}|^{ p-1}
|g^{k (\varepsilon)}_{s}|  \,dx\big)^{2}\leq
\big(
\int_{\bR^{d}}|u^{(\varepsilon)}_{s}|^{ p-1}
|g^{ (\varepsilon)}_{s}|_{\ell_{2}}\,dx\big)^{2}.
\end{equation}
By  H\"older's inequality
the right-hand side of \eqref{4.22.1} is less than 
$$
\big(\int_{\bR^{d}}|u^{(\varepsilon)}_{s}|^{ p }\,
dx\big)^{2(p-1)/p}
\big(\int_{\bR^{d}}|g^{  (\varepsilon)}_{s}|^{p}\,dx
\big)^{2/p}
$$
$$
\leq\big(\int_{\bR^{d}}|g^{  (\varepsilon)}_{s}|^{p}\,dx
\big)^{2/p}\big(
\int_{\bR^{d}}\sup_{s\leq t }
|u^{(\varepsilon)}_{s}|^{ p }\,
dx\big)^{2(p-1)/p}.
$$
Here
$$
|u^{(\varepsilon)}_{s}|^{ p }\leq (|u_{s}|^{p})^{(\varepsilon)},
\quad 
\sup_{s\leq t }
|u^{(\varepsilon)}_{s}|^{ p }
\leq (\sup_{s\leq t }|u_{s}|^{p})^{(\varepsilon)},
$$
$$
\int_{\bR^{d}}(\sup_{s\leq t }|u_{s}|^{p})^{(\varepsilon)}\,dx
=\int_{\bR^{d}} \sup_{s\leq t }|u_{s}|^{p} \,dx,
$$

It follows from here and \eqref{5.9.7}
that the process
$$
\xi^{\varepsilon}_{t}:=
\int_{0}^{t}\sum_{k=1}^{\infty}
\big(\int_{\bR^{d}}|u^{(\varepsilon)}_{s}|^{ p-1}
|g^{k (\varepsilon)}_{s}|_{\ell_{2}}\,dx\big)^{2} \,ds 
$$
is well defined, $\cF_{t}$ adapted,
 and is continuous in $t$ (a.s.). Hence,
the stopping times
$$
\tau_{n}=\tau\wedge\inf\{t\geq0:
\xi^{\varepsilon}_{t}\geq n\},
$$
$n=1,2,...$, are well defined, $\tau_{n}\uparrow\infty$,
and, obviously, \eqref{4.16.3} holds.

Estimates \eqref{4.18.2} show that there is no trouble in applying
the deterministic Fubini's theorem to the integrals with respect
to $ds$ in \eqref{4.12.1}.
Estimate \eqref{4.16.3} implies
 that for each fixed $t$ we can apply the stochastic
Fubini's theorem  to the stochastic
term in \eqref{4.12.1} with $t\wedge\tau_{n}$
in place of $t$. Hence we obtain that with probability one
$$
\int_{\bR^{d}}|u^{(\varepsilon)}_{t\wedge\tau_{n}}|^{p}\,dx
=\int_{\bR^{d}}|u^{(\varepsilon)}_{0}|^{p}\,dx
+\int_{0}^{t\wedge\tau_{n}}
 \big[p\int_{\bR^{d}}|u^{(\varepsilon)}_{s}|^{p-2}
u^{(\varepsilon)}_{s}f^{(\varepsilon)}_{s}\,dx\,ds
$$
$$
+(1/2)p(p-1)\int_{\bR^{d}}|u^{(\varepsilon)}_{s}|^{p-2}
|g^{(\varepsilon)}_{s}|_{\ell_{2}}^{2}\,dx
\big]\,ds
$$
$$
+p
\int_{0}^{t\wedge\tau_{n}}\int_{\bR^{d}}|u^{(\varepsilon)}_{s}|^{p-2}
u^{(\varepsilon)}_{s}
g^{k(\varepsilon)}_{s}\,dx\,dw^{k}_{s}.
$$
Since $\tau_{n}\uparrow\infty$, we have also
$$
\int_{\bR^{d}}|u^{(\varepsilon)}_{t }|^{p}\,dx
=\int_{\bR^{d}}|u^{(\varepsilon)}_{0}|^{p}\,dx
+\int_{0}^{t }
 \big[p\int_{\bR^{d}}|u^{(\varepsilon)}_{s}|^{p-2}
u^{(\varepsilon)}_{s}f^{(\varepsilon)}_{s}\,dx\,ds
$$
$$
+(1/2)p(p-1)\int_{\bR^{d}}|u^{(\varepsilon)}_{s}|^{p-2}
|g^{(\varepsilon)}_{s}|_{\ell_{2}}^{2}\,dx
\big]\,ds
$$
\begin{equation}
                                                     \label{4.18.1}
+p
\int_{0}^{t }\int_{\bR^{d}}|u^{(\varepsilon)}_{s}|^{p-2}
u^{(\varepsilon)}_{s}
g^{k(\varepsilon)}_{s}\,dx\,dw^{k}_{s}
\end{equation}
(a.s.) for each $t$.

We now pass to the limit as $\varepsilon\to0$
in \eqref{4.18.1}. Observe that for any $h\in L_{p}$
we have $\|h^{(\varepsilon)}\|_{L_{p}}
\leq\|h\|_{L_{p}}$ and $h^{(\varepsilon)}\to h$ in $L_{p}$
as $\varepsilon\to0$. Therefore, the left-hand side of 
\eqref{4.18.1} tends to the left-hand side of 
\eqref{4.18.5} with ($\tau=\infty$) for all $(\omega,t)$.
The same is true (a.s.) for the first term on the right
in \eqref{4.18.1}.

To prove the  convergence in probability of stochastic integrals
it suffices to prove that 
\begin{equation}
                                                         \label{4.19.2}
\int_{0}^{t}\sum_{k=1}^{\infty}
\big(\int_{\bR^{d}}\big(|u^{(\varepsilon)}_{s}|^{p-2}
u^{(\varepsilon)}_{s}g^{k(\varepsilon)}_{s}-
|u _{s}|^{p-2}u _{s}g^{k }_{s}\big)\,dx\big)^{2}ds\to0
\end{equation}
as $\varepsilon\to0$ (a.s.). Notice that 
for $s\leq t$
by Minkowski's and H\"older's
inequalities
$$
 \sum_{k=1}^{\infty}
\big(\int_{\bR^{d}}
|u _{s}|^{p-2}u _{s} (g^{k(\varepsilon)}_{s}-
 g^{k }_{s}\big)\,dx\big)^{2}ds\leq
\big(\int_{\bR^{d}}
|u _{s}|^{p-1}|g^{ (\varepsilon)}_{s}-
 g  _{s}|_{\ell_{2}}\,dx\big)^{2}
$$
$$
\leq \sup_{r\leq t}\|u_{r}\|^{2(p-1)}_{L_{p}}
\|g^{ (\varepsilon)}_{s}-
 g  _{s}\|^{2  }_{L_{p}}.
$$
The integral of the last term in $s$ tends to zero as $\varepsilon
\to 0$ (a.s.) owing to the above mentioned
 properties of
mollifiers and the fact that $\sup_{r\leq t}\|u_{r}\|_{L_{p}}<
\infty$ (for all $\omega$).

Hence to prove \eqref{4.19.2} it suffices to show that
$$
J^{\varepsilon}:=\int_{0}^{t}\sum_{k=1}^{\infty}
\big(\int_{\bR^{d}} (|u^{(\varepsilon)}_{s}|^{p-2}
u^{(\varepsilon)}_{s} -
|u _{s}|^{p-2}u _{s})g^{k (\varepsilon)}_{s}\,dx\big)^{2}ds 
$$
tends to zero (a.s.). By Minkowski's inequality 
\begin{equation}
                                                         \label{4.19.4}
J^{\varepsilon}\leq\int_{0}^{t} (I_{s}^{\varepsilon})^{2}\,ds,
\end{equation}
where
$$
I_{s}^{\varepsilon}:=
\int_{\bR^{d}}
|\,|u^{(\varepsilon)}_{s}|^{p-2}
u^{(\varepsilon)}_{s} -
|u _{s}|^{p-2}u _{s}|\,
|g^{  (\varepsilon)}_{s}|_{\ell_{2}}\,dx.
$$
Observe that on a set of full probability for almost any $s$
we have
$$
|g^{  (\varepsilon)}_{s}-g_{s}|_{\ell_{2}}\to0,\quad
|g^{  (\varepsilon)}_{s}|_{\ell_{2}}\to|g _{s}|_{\ell_{2}},
\quad
u^{  (\varepsilon)}_{s}\to u_{s}
$$
in $L_{p}$. In particular,
$$
\|\,|u^{(\varepsilon)}_{s}|^{p-2}
u^{(\varepsilon)}_{s}\|_{L_{p/(p-1)}}=\|u^{(\varepsilon)}_{s}\|_{L_{p}}\to
\|\,|u _{s}|^{p-2}
u _{s}\|_{L_{p/(p-1)}}.
$$
By Lemma \ref{lemma 1.15.1}
$$
|u^{(\varepsilon)}_{s}|^{p-2}
u^{(\varepsilon)}_{s}-|u _{s}|^{p-2}
u _{s}\to0
$$
in $L_{p/(p-1)}$. It follows that
 $
I_{s}^{\varepsilon} \to 0
 $
as $\varepsilon\to0$ for almost all $s$ if
$\omega\in\Omega'$ with $P(\Omega')=1$. Furthermore,
$$
I_{s}^{\varepsilon}\leq \sup_{r\leq t}
\|u^{ (\varepsilon)}_{r}\|^{  p-1 }_{L_{p}}
\|g^{ (\varepsilon)}_{s} \| _{L_{p}}
+\sup_{r\leq t}
\|u _{r}\|^{ p-1 }_{L_{p}}
\|g^{ (\varepsilon)}_{s} \| _{L_{p}} 
$$
$$
\leq 2\sup_{r\leq t}
\|u _{r}\|^{ p-1 }_{L_{p}}
\|g _{s} \| _{L_{p}},
$$
which is square integrable over $[0,t]$
(a.s.).
By the dominated convergence theorem
and \eqref{4.19.4} we have $J^{\varepsilon}\to0$ (a.s.) 
yielding the desired convergence in probability of the 
stochastic integrals in \eqref{4.18.1} as $\varepsilon\to0$.

The integral with respect to $ds$ on the right in \eqref{4.18.1}
presents no
difficulty owing to  Corollary \ref{corollary 1.15.1}
and is treated similarly to what is done
above after~\eqref{4.19.4}.

Thus, for each $t$ equation \eqref{4.18.5} holds
with probability one. Since both parts
are continuous in $t$, it also holds for all
$t$ at once on the set of full probability and
 this finally
brings the proof of the lemma to an end.

\mysection{Proof of Theorem \protect\ref{theorem 12.3.1}
and Corollary \protect\ref{corollary 4.19.1}}
 
                                              \label{section 5.13.2}

First we prove the theorem.  We use the notation $h^{(\varepsilon)}$
as in the proof of Lemma \ref{lemma 4.17.1}
taking there a nonnegative $\zeta\in C^{\infty}_{0}$
 with unit integral. By substituting $\phi*\bar{\zeta}_{\varepsilon}$,
where $\bar{\zeta}_{\varepsilon}(x)=\zeta _{\varepsilon}(-x)$,
 in place of $\phi$
in \eqref{1.2.1} we see that 
$u^{(\varepsilon)}_{t}$ 
satisfies 
$$
(u_{t\wedge\tau}^{(\varepsilon)},\varphi)=(u_{0}^{(\varepsilon)},\varphi)
+\sum_{k=1}^{\infty}\int_{0}^{t}I_{s\leq\tau}
(g^{k(\varepsilon)}_{s},\varphi)\,dw^{k}_{s}  
$$
\begin{equation}
                                                     \label{5.12.1}
+\int_{0}^{t}I_{s\leq\tau}\big(-(f^{i(\varepsilon)}_{s},D_{i}\varphi)
+(f^{0(\varepsilon)}_{s},\varphi)\big)\,ds,
\end{equation}
which is
\eqref{4.18.6} with
$u_{0}^{(\varepsilon)}$ and $g^{k(\varepsilon)}_{t}$
in place of $u_{0}$ and $g^{k}_{t}$, respectively, and with
$$
 D_{i}(f^{i})^{(\varepsilon)}+(f^{0})^{(\varepsilon)}
$$
in place of $f_{t}$.

From   Lemma \ref{lemma 4.17.1}
 we obtain that, for an $\Omega^{\varepsilon}$
with $P(\Omega^{\varepsilon})=1$, $u^{(\varepsilon)}_{t\wedge\tau}
I_{\Omega^{\varepsilon}}$ is a continuous
$L_{p}$-valued $\cF_{t}$-adapted process on $[0,\infty)$  
and the corresponding counterpart of
\eqref{4.18.5} holds, integrating by parts in which leads to
the fact that with probability one for all $t\geq0$
$$
\int_{\bR^{d}}|u^{(\varepsilon)}_{t\wedge\tau}|^{p}\,dx
=\int_{\bR^{d}}|u^{(\varepsilon)}_{0}|^{p}\,dx
+p
\int_{0}^{t\wedge\tau }\int_{\bR^{d}}|u^{(\varepsilon)}_{s}|^{p-2}
u^{(\varepsilon)}_{s}
g^{k (\varepsilon)}_{s}\,dx\,dw^{k}_{s}
$$
$$+
\int_{0}^{t\wedge\tau }
\big( \int_{\bR^{d}}\big[p|u^{(\varepsilon)}_{s}|^{p-2}u^{(\varepsilon)}_{s}
f^{0(\varepsilon)}_{s}
-p(p-1)|u^{(\varepsilon)}_{s}|^{p-2}
f^{i(\varepsilon)}_{s}D_{i}u^{(\varepsilon)}_{s}
$$
\begin{equation}
                                            \label{4.19.05}
+(1/2)p(p-1)|u^{(\varepsilon)}_{s}|^{p-2}
|g^{(\varepsilon)}_{s}|_{\ell_{2}}^{2}
\big]\,dx\big)\,ds.
\end{equation}

We take the supremums with respect
to $t$ of both parts and
 repeat a standard  argument which was introduced by E. Pardoux.
We will be using \eqref{4.18.2}  
and the fact that,
by the inequality $a^{p-2}bc\leq a^{p}+b^{p}+c^{p}$, $a,b,c\geq0$,
we have
$$
\int_{\bR^{d}}|u^{(\varepsilon)}_{s}|^{p-2}
f^{i(\varepsilon)}_{s}D_{i}u^{(\varepsilon)}_{s}\,dx
\leq
\frac{\gamma^{p/(p-2)}}{T}\int_{\bR^{d}}|u^{(\varepsilon)}_{s}|^{p }\,dx
$$
\begin{equation}
                                                 \label{5.13.2}
+\frac{T^{(p-2)/2}}{\gamma^{p/2}}
\int_{\bR^{d}}[|f^{(\varepsilon)}_{s}|^{p}
+|Du^{(\varepsilon)}_{s}|^{p}]\,dx,
\end{equation}
where $f^{(\varepsilon)}=(f^{1( \varepsilon)},...,f^{d( \varepsilon)})$
and $\gamma>0$ is any number. We also use \eqref{4.22.1}
and  the  Burkholder-Davis-Gundy
inequalities.
Then for an appropriate choice of the parameter $\gamma$ 
we find from \eqref{4.19.05} that
$$
E\sup_{t\leq\tau\wedge T}\|u^{(\varepsilon)}_{t }\|_{L_{p}}^{p}
\leq E\|u^{(\varepsilon)}_{0}\|_{L_{p}}^{p}+
(1/4)E\sup_{t\leq\tau\wedge T}\|u^{(\varepsilon)}_{t }\|_{L_{p}}^{p}
$$
$$
+NT^{(p-2)/2}E\int_{0}^{T\wedge\tau}(\|g^{(\varepsilon)}_{s}\|_{L_{p}}^{p }
+\|f^{(\varepsilon)}_{s}\|_{L_{p}}^{p } 
+\|Du^{(\varepsilon)}_{s}\|_{L_{p}}^{p })\,ds
$$
$$
+NT^{p-1}E 
\int_{0}^{T\wedge\tau}\|f^{0(\varepsilon)}_{s}\|_{L_{p}}^{p }\,ds
+NE\big(\int_{0}^{T\wedge\tau}
\big(\int_{\bR^{d}}|u^{(\varepsilon)}_{s} |^{p-1}
 |g^{(\varepsilon)}_{s}|_{\ell_{2}} \,dx\big)^{2}\,ds\big)^{1/2},
$$
where the last expectation is estimated by
$$
E\sup_{t\leq\tau\wedge T}\|u^{(\varepsilon)}_{t }\|_{L_{p}}^{p-1}
\big(\int_{0}^{T\wedge\tau} 
\|g^{(\varepsilon)}_{s}\|_{L_{p}}^{2}\,ds\big)^{1/2}
$$
$$
\leq T^{(p-2)/(2p)}
E\sup_{t\leq\tau\wedge T}\|u^{(\varepsilon)}_{t }\|_{L_{p}}^{p-1}
\big(\int_{0}^{T\wedge\tau} 
\|g^{(\varepsilon)}_{s}\|_{L_{p}}^{p}\,ds\big)^{1/p}
$$
$$
\leq(1/4)E\sup_{t\leq\tau\wedge T}\|u^{(\varepsilon)}_{t }\|_{L_{p}}^{p}
+NT^{(p-2)/2}E
\int_{0}^{T\wedge\tau} 
\|g^{(\varepsilon)}_{s}\|_{L_{p}}^{p}\,ds.
$$
Hence,
$$
E\sup_{t\leq\tau\wedge T}\|u^{(\varepsilon)}_{t }\|_{L_{p}}^{p}
\leq E\|u^{(\varepsilon)}_{0}\|_{L_{p}}^{p}+
(1/2)E\sup_{t\leq\tau\wedge T}\|u^{(\varepsilon)}_{t }\|_{L_{p}}^{p}
$$
$$
+NT^{(p-2)/2}E\int_{0}^{T\wedge\tau}(\|g^{(\varepsilon)}_{s}\|_{L_{p}}^{p }
+\|f^{(\varepsilon)}_{s}\|_{L_{p}}^{p } 
+\|Du^{(\varepsilon)}_{s}\|_{L_{p}}^{p })\,ds
$$
$$
+NT^{p-1}E 
\int_{0}^{T\wedge\tau}\|f^{0(\varepsilon)}_{s}\|_{L_{p}}^{p }\,ds.
$$
Upon collecting like terms we come to
$$
E\sup_{t\leq\tau\wedge T}\|u^{(\varepsilon)}_{t }\|_{L_{p}}^{p}
\leq 2E\|u^{(\varepsilon)}_{0}\|_{L_{p}}^{p} 
+NT^{p-1}E 
\int_{0}^{T\wedge\tau}\|f^{0(\varepsilon)}_{s}\|_{L_{p}}^{p }\,ds   
$$
\begin{equation}
                                                     \label{4.16.2} 
+NT^{(p-2)/2}E\int_{0}^{T\wedge\tau}(\|g^{(\varepsilon)}_{s}\|_{L_{p}}^{p }
+\|f^{(\varepsilon)}_{s}\|_{L_{p}}^{p } 
+\|Du^{(\varepsilon)}_{s}\|_{L_{p}}^{p })\,ds.
\end{equation}
One can lawfully object that the last step leads to estimate
\eqref{4.16.2} only if its left-hand side is finite.
However, by Lemma \ref{lemma 4.17.1} the process 
$\|u^{(\varepsilon)}_{t }\|_{L_{p}}^{p}$ is a continuous $\cF_{t}$-adapted
process 
which starts at $\|u^{(\varepsilon)}_{0}\|_{L_{p}}^{p}$
and we can stop it 
at time $\tau_{n}$
when it first reaches the level 
$\|u^{(\varepsilon)}_{0}\|_{L_{p}}^{p} +n$ with $n>0$ or at time 
  $\tau$ whichever comes first. Then 
$$
E\sup_{t\leq\tau_{n}\wedge T}\|u^{(\varepsilon)}_{t }\|_{L_{p}}^{p}
\leq E\|u^{(\varepsilon)}_{0}\|_{L_{p}}^{p} +n<\infty.
$$
Hence,
the left-hand side
of \eqref{4.16.2} will be finite if we replace there
$\tau$ with $\tau_{n}$. Therefore, thus modified
\eqref{4.16.2}   holds and sending $n$ to infinity
yields \eqref{4.16.2} as is.

By applying this result to $u^{(\varepsilon_{1})}_{t}
-u^{(\varepsilon_{2})}_{t}$ we conclude that 
$$
E\sup_{t\leq\tau\wedge T}
\|u^{(\varepsilon_{1})}_{t}
-u^{(\varepsilon_{2})}_{t}\|^{p}_{L_{p}}\to0
$$
as $\varepsilon_{1},\varepsilon_{2}\to0$. It follows that
there exists a function $v_{t}=v_{t}(\omega,x)$, $0\leq t\leq\tau(\omega),
t<\infty,x\in\bR^{d}$, such that
  $v_{t\wedge\tau}$   is   continuous in $t$ and $\cF_{t}$-adapted
as an $L_{p}$-valued function and
\begin{equation}
                                         \label{4.20.2}
E\sup_{t\leq\tau\wedge T}
\|u^{(\varepsilon )}_{t}
-v_{t}\|^{p}_{L_{p}}\to0
\end{equation}
as $\varepsilon\to0$. 
In particular, in probability, 
for any $\phi\in C^{\infty}_{0}$,   we have
\begin{equation}
                                                  \label{4.18.7}
(u^{(\varepsilon )}_{t\wedge\tau},\phi)\to (v_{t\wedge\tau},\phi)
\end{equation}
uniformly on $[0,T]$ for any $T\in[0,\infty)$. Also in probability,
$$
\big|\int_{0}^{T} (|(f^{i(\varepsilon)}_{s}-
f^{i }_{s},D_{i}\varphi)|
+|(f^{0(\varepsilon)}_{s}-f^{0}_{s},\varphi)|)\,ds\big|^{p}
$$
$$
\leq N \int_{0}^{T}\sum_{j=0}^{d}\| f^{j(\varepsilon)}_{s}-
f^{i }_{s}\|_{L_{p}}^{p}\,ds\to0,
$$
and (cf. \eqref{5.12.2})
$$
\big|\int_{0}^{T}\sum_{k=1}^{\infty} (g^{k(\varepsilon)}_{s}
-g^{k}_{s},\varphi) ^{2}\,ds\big|^{p/2}
\leq N\int_{0}^{T}\|g^{(\varepsilon)}_{s}-g_{s}\|^{p}_{L_{p}} \,ds
\to 0.
$$
Therefore, we can pass to the limit in \eqref{5.12.1}
and conclude that \eqref{12.3.1} holds with $v$ in place of $u$.
 The same  argument as in the proof of Lemma
\ref{lemma 4.17.1} now shows that with probability one
the generalized functions $v_{t\wedge\tau}$ and $u_{t\wedge\tau}$ coincide
for all $t\in[0,\infty)$. 
This proves the assertion (i) of the theorem.
 After that \eqref{4.11.5} 
is obtained by sending $\varepsilon\to0$ 
in \eqref{4.16.2}. Finally,
the argument in  the proof of Lemma
\ref{lemma 4.17.1}  can also be repeated almost
 literally to obtain formula \eqref{4.19.5} from \eqref{4.19.05}.
 The theorem is proved.

{\bf Proof of Corollary \ref{corollary 4.19.1}}.
Denote by $J(\tau)$ the left-hand side of \eqref{12.3.2}.
Estimates similar to \eqref{4.18.2} and \eqref{5.13.2} show that
$J(\tau)<\infty$ and if we have a sequence of stopping times
$\tau_{n}\uparrow\tau$, then $J(\tau_{n})\to J(\tau)$ as $n
\to\infty$. By taking a sequence which localizes the stochastic
integral in \eqref{4.19.5}, then taking
expectations of both sides of \eqref{4.19.5},  and,
finally, using Fatou's lemma,
we obtain the first assertion of the corollary. If $\tau$ is bounded,
then the described procedure will yield the second assertion
as well, due to \eqref{4.11.5} and the dominated convergence theorem.

\end{document}